\documentstyle[12pt]{article}

\catcode`\@=11
\begin{document}
\centerline{\large\bf Super Yangian Double $DY(gl(1|1))$ With Center Extension }
\vspace{0.8cm}
\centerline{\sf $^a$Jin-fang Cai, $^a$Guo-xing Ju $^{bc}$S.K. Wang and  
$^a$Ke Wu } 
\baselineskip=13pt
\vspace{0.5cm}
\centerline{$^a$ Institute of Theoretical Physics, Academia Sinica, }
\baselineskip=12pt
\centerline{ Beijing, 100080, P. R. China }
\vspace{0.3cm}
\baselineskip=13pt
\centerline{ $^b$ CCAST (World Laboratory), P.O. Box 3730, Beijing, 100080, 
  P. R. China  }
\baselineskip=12pt
\vspace{0.3cm}
\centerline{ $^c$ Institute of Applied Mathematics, Academia Sinica,}
\baselineskip=12pt
\centerline{Beijing, 100080, P. R. China }
\vspace{0.9cm}
\begin{abstract} {We use the super version of RS construction method and 
Gauss decomposition to obtain Drinfel$^{\prime}$d's  currents realization of 
super Yangian double $DY(gl(1|1))$ with center extension}
\end{abstract}
\baselineskip=14pt
\vspace{1cm}
Yangians and quantum affine algebras are quantum algebras which are related 
respectively with rational and  trigonometric solutions of Yang-Baxter equation
\cite{DRI1}.
The explicit isomorphism between two realizations of quantum affine 
algebras and Yangians  given by Drinfel$^{\prime}$d \cite{DRI2} and RS \cite{RS}
was established by Ding and Frenkel \cite{DF} using Gauss decomposition. 
In this paper, we use super version of RS construction method to define 
super Yangian double $DY(gl(1|1))$ \cite{C1} with center extension, and get its
 Drinfeld's currents realization by Guass decomposition.

By $V$ we denote a ${\bf Z}_2$-graded two-dimension vector space 
(graded auxiliary space), and set the first and second basis of $V$ is 
even and odd respectively. In graded case, we must use this 
tensor product form: $(A \otimes B)(C \otimes D) =(-1)^{P(B)P(C)}
AC \otimes BD $, here $P(A)=0,1$ when $A$ is bosonic and Fermionic 
respectively. Then the graded(super) Yang-Baxter equation(YBE) takes
 this form\cite{LS}: 
\begin{eqnarray}
&&\eta _{12}R_{12}(u) \eta _{13}R_{13}(u+v) \eta _{23}R_{23}(v) 
= \eta _{23}R_{23}(v) \eta _{13}R_{13}(u+v) \eta _{12}R_{12}(u),
\end{eqnarray}
where $R(u) \in End(V\otimes V) $ and it must obey weight conservation 
condition: $R_{ij,kl}\neq 0 $ only when $i+j = k+l $.
 $\eta _{ik,jl}=(-1)^{(i-1)(k-1)} \delta _{ij} \delta _{lk}$ .
It's easily proved that the following $R(u)$ matrix is a rational solution
of super YBE:

\begin{equation}
R(u)=\frac{1}{u+\hbar}(uI+{\cal P})=\left(\begin{array}{cccc}
1&0&0&0\\
0&\frac{u}{u+\hbar}&\frac{\hbar}{u+\hbar}&0\\
0&\frac{\hbar}{u+\hbar}&\frac{u}{u+\hbar}&0\\
0&0&0&\frac{u-\hbar}{u+\hbar} \end{array}\right)
\end{equation}

This sulution of graded YBE satisfied unitarity condition: $R(u)R_{21}(-u)=I$.
From this rational solution  of graded YBE, we can define the 
super Yangian double $DY(gl(1|1))$ with a central extension  employing 
super version of RS method \cite{RS}.  Super Yangian Double   \\
$DY\left( gl(1|1)\right)$ is a Hopf algebra generated 
by $ \{ l_{ij}^k \vert 1\leq i,j \leq 2, k\in {\bf Z }\} $ which obey the 
following relations:
\begin{eqnarray}
&&R(u-v)L_1^{\pm}(u)\eta L_2^{\pm}(v)\eta=\eta L_2^{\pm}(v)
  \eta L_1^{\pm}(u)R(u-v)  \label{rll1} \\
&&R(u_--v_+)L_1^+(u)\eta L_2^-(v)\eta=\eta L_2^-(v)
  \eta L_1^+(u)R(u_+-v_-)  \label{rll2}
\end{eqnarray}
here $u_{\pm}=u\pm \frac{\hbar}{4}c$ and 
we have used standard notation: $L_1^{\pm}(u)=L^{\pm}(u)\otimes {\bf 1},
 L_2^{\pm}(u)={\bf 1} \otimes
L^{\pm}(u) $ and $L^{\pm}(u)=\left( l_{ij}^{\pm}(u) \right)
 _{1 \leq i,j \leq 2}$,
  $l_{ij}^{\pm}(u)$ are generate functions of $l_{ij}^k$ :
\begin{equation}
 l_{ij}^{+}(u)=\delta _{ij}-\hbar\sum_{k\geq 0}l_{ij}^k u^{-k-1},
  {\hspace{0.5cm}}
 l_{ij}^{-}(u)=\delta _{ij}+\hbar\sum_{k < 0}l_{ij}^k u^{-k-1}.
\end{equation}
 Hopf structure of $DY\left(gl(1\vert 1)\right)$  are defined as follows:
\begin{eqnarray}
&&\triangle \left(l_{ij}^{\pm}(u)\right) =\sum_{k=1,2} l_{kj}^{\pm}(u \pm \frac{\hbar}{4}c_2)
 \otimes l_{ik}^{\pm}(u\mp \frac{\hbar}{4}c_1)(-1)^{(k+i)(k+j)}, \label{cop} \\
&& \epsilon\left(l_{ij}^{\pm}(u)\right)=\delta _{ij} ,{\hspace{1.cm}}
S\left( ^{st} L^{\pm}(u) \right) =\left[ ^{st} L^{\pm}(u) \right]^{-1}.
\label{es}
\end{eqnarray} 
here $\left[^{st} L^{\pm}(u) \right]_{ij} =(-1)^{i+j}l_{ji}^{\pm}(u) $ .

From Ding-Frenkel theorem\cite{DF} on two realizations of Yangians  and quantum 
affine algebras, $L^{\pm}(u)$ have the following decomposition:
\begin{eqnarray}
L^{\pm}(u)&=&\left(\begin{array}{cc}1&0\\f^{\pm}(u)&1\end{array}\right)
\left(\begin{array}{cc}k_1^{\pm}(u)&0\\0&k_2^{\pm}(u)\end{array}\right)
\left(\begin{array}{cc}1&e^{\pm}(u)\\0&1\end{array}\right)  \nonumber\\
&=&\left(\begin{array}{cc}k_1^{\pm}(u)&k_1^{\pm}(u)e^{\pm}(u)\\f^{\pm}(u)
k_1^{\pm}(u)&k_2^{\pm}(u)+f^{\pm}(u)k_1^{\pm}(u)e^{\pm}(u)\end{array}\right) 
\end{eqnarray}
and the  inversions of $L^{\pm}(u)$ can be writed as:
\begin{eqnarray}
L^{\pm}(u)^{-1}&=&\left(\begin{array}{cc}1&-e^{\pm}(u)\\0&1\end{array}\right)
\left(\begin{array}{cc}k_1^{\pm}(u)^{-1}&0\\0&k_2^{\pm}(u)^{-1}\end{array}\right)
\left(\begin{array}{cc}1&0\\-f^{\pm}(u)&1\end{array}\right)  \nonumber\\
&=&\left(\begin{array}{cc}k_1^{\pm}(u)^{-1}+
e^{\pm}(u)k_2^{\pm}(u)^{-1}f^{\pm}(u)&-e^{\pm}(u)k_2^{\pm}(u)^{-1}\\
-k_2^{\pm}(u)^{-1}f^{\pm}(u)&k_2^{\pm}(u)^{-1}\end{array}\right) 
\end{eqnarray}
Let
\begin{equation}
X^+(v)=e^+(v_-)-e^-(v_+)  {\hspace{1cm}} X^-(v)=f^+(v_+)-f^-(v_-)
\end{equation}

To calculate the (anti-)commutation relations of $X^{\pm}(u)$ and 
$k_i^{\pm}(v)$ (i=1,2), we will use the following equivalent forms of 
eq.(\ref{rll1})(\ref{rll2}):
\begin{eqnarray}
L_1^{\pm}(u)^{-1}\eta L_2^{\pm}(v)^{-1}\eta R(u-v)&=
&R(u-v)\eta L_2^{\pm}(v)^{-1} \eta L_1^{\pm}(u)^{-1} \label{rll3}  \\
L_1^+(u)^{-1}\eta L_2^-(v)^{-1}\eta R(u_--v_+)&=&R(u_+-v_-)\eta L_2^-(v)^{-1}
 \eta L_1^+(u)^{-1} \label{rll4} \\
L_1^{\pm}(u)R(u-v)\eta L_2^{\pm}(v)^{-1} \eta &=&\eta L_2^{\pm}(v)^{-1} 
\eta R(u-v)L_1^{\pm}(u)\label{rll5} \\
L_1^+(u)R(u_+-v_-)\eta L_2^-(v)^{-1} \eta &=&\eta L_2^-(v)^{-1} 
\eta R(u_--v_+)L_1^+(u) \label{rll6}\\
L_1^{\pm}(u)^{-1}R_{21}(v-u)\eta L_2^{\pm}(v) \eta&=
&\eta L_2^{\pm}(v) \eta R_{21}(v-u)L_1^{\pm}(u)^{-1}  \label{rll7}\\
L_1^+(u)^{-1}R_{21}(v_+-u_-)\eta L_2^-(v) \eta&=
&\eta L_2^-(v) \eta R_{21}(v_--u_+)L_1^+(u)^{-1} \label{rll8}
\end{eqnarray}

From (\ref{rll1}) (\ref{rll2}) and (\ref{rll3})---(\ref{rll8}), 
we can get all relations between 
$k_1^{\pm}(u)$ and $  k_2^{\pm}(v) $ :
\begin{eqnarray}
&&[k_1^{\pm}(u) ~,~k_1^{\pm}(v)]=[k_1^+(u) ~,~k_1^-(v)]=0 \label{k1}\\
&&[k_1^{\pm}(u) ~,~k_2^{\pm}(v)]=[k_2^{\pm}(u) ~,~k_2^{\pm}(v)]=0  \label{k2}\\
&&\frac{u_{\pm}-v_{\mp}}{u_{\pm}-v_{\mp}+\hbar}k_1^{\pm}(u)k_2^{\mp}(v)^{-1}
=\frac{u_{\mp}-v_{\pm}}{u_{\mp}-v_{\pm}+\hbar}k_2^{\mp}(v)^{-1}
k_1^{\pm}(u)  \label{k3}\\
&&\frac{u_--v_+-\hbar}{u_--v_++\hbar}k_2^+(u)^{-1}k_2^-(v)^{-1}
=\frac{u_+-v_--\hbar}{u_+-v_-+\hbar}k_2^-(v)^{-1}k_2^+(u)^{-1} \label{k4}
\end{eqnarray}

From (\ref{rll1})(\ref{rll2}) and unitarity of $R$-matrix, 
we have the following relations between $k_1^{\pm}(u)$ 
and $e^{\pm}(v), f^{\pm}(v)$:

\begin{eqnarray}
&&k_1^{\pm}(u)k_1^{\pm}(v)e^{\pm}(v)-\frac{u-v}{u-v+\hbar}k_1^{\pm}(v)
e^{\pm}(v)k_1^{\pm}(u) \nonumber\\   
&&{\hspace{5cm}}-\frac{\hbar}{u-v+\hbar}k_1^{\pm}(v)k_1^{\pm}(u)
e^{\pm}(u)=0  \\
&&k_1^{\pm}(u)k_1^{\mp}(v)e^{\mp}(v)-\frac{u_{\pm}-v_{\mp}}{u_{\pm}-v_{\mp}
+\hbar}k_1^{\mp}(v)e^{\mp}(v)k_1^{\pm}(u) \nonumber\\   
&&{\hspace{5cm}}-\frac{\hbar}{u_{\pm}-v_{\mp}+\hbar}
k_1^{\mp}(v)k_1^{\pm}(u)e^{\pm}(u)=0  \\
&&f^{\pm}(v)k_1^{\pm}(v)k_1^{\pm}(u)-\frac{u-v}{u-v+\hbar}k_1(u)f^{\pm}(v)
k_1^{\pm}(v)  \nonumber\\   &&{\hspace{5cm}}
-\frac{\hbar}{u-v+\hbar}f^{\pm}(u)k_1^{\pm}(u)k_1^{\pm}(v)=0   \label{kf1}  \\
&&f^{\mp}(v)k_1^{\mp}(v)k_1^{\pm}(u)-\frac{u_{\mp}-v_{\pm}}{u_{\mp}-v_{\pm}+
\hbar}k_1^{\pm}(u)f^{\mp}(v)k_1^{\mp}(v)
 \nonumber\\   &&{\hspace{5cm}}
-\frac{\hbar}{u_{\mp}-v_{\pm}+\hbar}
f^{\pm}(u)k_1^{\pm}(u)k_1^{\mp}(v)=0   \label{kf2}
\end{eqnarray}
thus
\begin{eqnarray}
&&(u_{\pm}-v+\hbar)e^{\pm}(v_{\mp})-(u_{\pm}-v)k_1^{\pm}(u)^{-1}e^{\pm}(v_{\mp})
k_1^{\pm}(u)-\hbar e^{\pm}(u)=0  \\
&&(u_{\pm}-v+\hbar)e^{\mp}(v_{\pm})-(u_{\pm}-v)k_1^{\pm}(u)^{-1}e^{\mp}(v_{\pm})
k_1^{\pm}(u)-\hbar e^{\pm}(u)=0  \\
&&(u_{\mp}-v+\hbar)f^{\pm}(v_{\pm})-(u_{\mp}-v)k_1^{\pm}(u)f^{\pm}(v_{\pm})
k_1^{\pm}(u)^{-1}-\hbar f^{\pm}(u)=0 \\
&&(u_{\mp}-v+\hbar)f^{\mp}(v_{\mp})-(u_{\mp}-v)k_1^{\pm}(u)f^{\mp}(v_{\mp})
k_1^{\pm}(u)^{-1}-\hbar f^{\pm}(u)=0 
\end{eqnarray}
so
\begin{eqnarray}
&&k_1^{\pm}(u)^{-1}X^+(v)k_1{\pm}(u)=\frac{u_{\pm}-v+\hbar}{u_{\pm}-v}X^+(v) \\
&&k_1^{\pm}(u)X^-(v)k_1^{\pm}(u)^{-1}=\frac{u_{\mp}-v+\hbar}{u_{\mp}-v}X^-(v)
\end{eqnarray}

Now, we derive relations between  $k_2^{\pm}(u)$ and $X^{\pm}(v)$. 
From(\ref{rll3})(\ref{rll4}) and unitarity of $R$-matirx, we have
\begin{eqnarray}
&&\frac{u-v-\hbar}{u-v+\hbar}e^{\pm}(u)k_2^{\pm}(u)^{-1}k_2^{\pm}(v)^{-1}-
\frac{u-v}{u-v+\hbar}k_2^{\pm}(v)^{-1}e^{\pm}(u)k_2^{\pm}(u)^{-1}
 \nonumber\\   &&{\hspace{5cm}}
+\frac{\hbar}{u-v+\hbar}e^{\pm}(v)k_2^{\pm}(v)^{-1}k_2^{\pm}(u)^{-1}=0  \\
&&\frac{u_{\mp}-v_{\pm}-\hbar}{u_{\mp}-v_{\pm}+\hbar}e^{\pm}(u)
k_2^{\pm}(u)^{-1}k_2^{\mp}(v)^{-1}-\frac{u_{\pm}-v_{\mp}}{u_{\pm}-
v_{\mp}+\hbar}k_2^{\mp}(v)^{-1}e^{\pm}(u)k_2^{\pm}(u)^{-1} \nonumber\\   &&{\hspace{5cm}}+\frac{\hbar}{u_{\pm}-v_{\mp}+\hbar}e^{\mp}(v)k_2^{\mp}(v)^{-1}
k_2^{\pm}(u)^{-1}=0  \\
&&\frac{u-v-\hbar}{u-v+\hbar}k_2^{\pm}(v)^{-1}k_2^{\pm}(u)^{-1}f^{\pm}(u)
-\frac{u-v}{u-v+\hbar}k_2^{\pm}(u)^{-1}f^{\pm}(u)k_2^{\pm}(v)^{-1}  \nonumber\\   &&{\hspace{5cm}}
+\frac{\hbar}{u-v+\hbar}k_2^{\pm}(u)^{-1}k_2^{\pm}(v)^{-1}f^{\pm}(v)=0  \\
&&\frac{u_{\pm}-v_{\mp}-\hbar}{u_{\pm}-v_{\mp}+\hbar}k_2^{\mp}(v)^{-1}
k_2^{\pm}(u)^{-1}f^{\pm}(u)-\frac{u_{\mp}-v_{\pm}}{u_{\mp}-v_{\pm}+\hbar}
k_2^{\pm}(u)^{-1}f^{\pm}(u)k_2^{\mp}(v)^{-1} \nonumber\\  
 &&{\hspace{5cm}}+\frac{\hbar}{u_{\mp}-v_{\pm} 
+\hbar}k_2^{\pm}(u)^{-1}k_2^{\mp}(v)^{-1}f^{\mp}(v)=0 
\end{eqnarray}
Using relations between $k_2^+(u)$ and $k_2^-(v)$ (\ref{k2})(\ref{k4}), we have
\begin{eqnarray}
&&(u_{\mp}-v-\hbar)e^{\pm}(u_{\mp})-(u_{\mp}-v)k_2^{\pm}(v)^{-1}e^{\pm}(u_{\mp})
k_2^{\pm}(v)+\hbar e^{\pm}(v)=0 \\
&&(u-v_{\mp}-\hbar)e^{\pm}(u_{\mp})-(u-v_{\mp})k_2^{\mp}(v)^{-1}e^{\pm}(u_{\mp})
k_2^{\mp}(v)+\hbar e^{\mp}(v)=0  \\
&&(u_{\pm}-v-\hbar)f^{\pm}(u_{\pm})-(u_{\pm}-v)k_2^{\pm}(v)f^{\pm}(u_{\pm})
k_2^{\pm}(v)^{-1}+\hbar f^{\pm}(v)=0 \\
&&(u-v_{\pm}-\hbar)f^{\pm}(u_{\pm})-(u-v_{\pm})k_2^{\mp}(v)f^{\pm}(u_{\pm})
k_2^{\mp}(v)^{-1}+\hbar f^{\mp}(v)=0
\end{eqnarray}
then
\begin{eqnarray}
&&k_2^{\pm}(u)^{-1}X^+(v)k_2^{\pm}(u)=\frac{u_{\pm}-v+\hbar}{u_{\pm}-v}X^+(v)  \\
&&k_2^{\pm}(u)X^-(v)k_2^{\pm}(u)^{-1}=\frac{u_{\mp}-v+\hbar}{u_{\mp}-v}X^-(v)
\end{eqnarray}

Now, we calculate the relations between  $X^{\pm}(u)$ 
and   $X^{\pm}(v)$ . From (\ref{rll1})(\ref{rll2}) and unitarity of
 $R$-matirx , 
we have
\begin{eqnarray}
&&k_1^{\pm}(u)e^{\pm}(u)k_1^{\pm}(v)e^{\pm}(v)+\frac{u-v-\hbar}
{u-v+\hbar}k_1^{\pm}(v)e^{\pm}(v)k_1^{\pm}(u)e^{\pm}(u)=0  \\
&&k_1^{\pm}(u)e^{\pm}(u)k_1^{\mp}(v)e^{\mp}(v)+\frac{u_{\pm}-v_{\mp}-\hbar}
{u_{\pm}-v_{\mp}+\hbar}k_1^{\mp}(v)e^{\mp}(v)k_1^{\pm}(u)e^{\pm}(u)=0  \\
&&\frac{u-v-\hbar}{u-v+\hbar}f^{\pm}(u)k_1^{\pm}(u)f^{\pm}(v)k_1^{\pm}(v)
+f^{\pm}(v)k_1^{\pm}(v)f^{\pm}(u)k_1^{\pm}(u)=0  \\
&&\frac{u_{\mp}-v_{\pm}-\hbar}{u_{\mp}-v_{\pm}+\hbar}f^{\pm}(u)k_1^{\pm}
(u)f^{\mp}(v)k_1^{\mp}(v)+f^{\mp}(v)k_1^{\mp}(v)f^{\pm}(u)k_1^{\pm}(u)=0 
\end{eqnarray}
Using the above relations between $X^{\pm}(u)$ (or $ e^{\pm}(u)$ and $f^{\pm}(u)$)
 and $k_1^{\pm}(v)$, it's easily 
to get the following relations;
\begin{eqnarray}
&&X^+(u)X^+(v)+X^+(v)X^+(u)=0  \\
&&X^-(u)X^-(v)+X^-(v)X^-(u)=0
\end{eqnarray}

From (\ref{rll1})(\ref{rll2}) and unitarity of $R$-matirx, we have following relations  :
\begin{eqnarray}
&&\frac{u-v}{u-v+\hbar}k_1^{\pm}(u)e^{\pm}(u)f^{\pm}(v)k_1^{\pm}(v)
+\frac{u-v}{u-v+\hbar}f^{\pm}(v)k_1^{\pm}(v)k_1^{\pm}(u)e^{\pm}(u) \nonumber  \\
&&{\hspace{4cm}}-\frac{\hbar}{u-v+\hbar}\left(k_2^{\pm}(u)+f^{\pm}(u)
k_1^{\pm}(u)e^{\pm}(u)
\right)k_1^{\pm}(v) \nonumber  \\
&&{\hspace{4cm}}-\frac{\hbar}{u-v+\hbar}\left(k_2^{\pm}(v)+f^{\pm}(v)k_1^{\pm}(v)
e^{\pm}(v)\right)k_1^{\pm}(u)=0  \\
&&\frac{u_{\mp}-v_{\pm}}{u_{\mp}-v_{\pm}+\hbar}k_1^{\pm}(u)e^{\pm}(u)
f^{\mp}(v)k_1^{\mp}(v) 
+\frac{u_{\pm}-v_{\mp}}{u_{\pm}-v_{\mp}+\hbar}f^{\mp}(v)k_1^{\mp}(v)
k_1^{\pm}(u)e^{\pm}(u)  \nonumber \\
&&{\hspace{4cm}}-\frac{\hbar}{u_{\mp}-v_{\pm}+\hbar}\left(k_2^{\pm}(u)
+f^{\pm}(u)k_1^{\pm}(u)e^{\pm}(u)\right)k_1^{\mp}(v) \nonumber  \\
&&{\hspace{4cm}}-\frac{\hbar}{u_{\pm}-v_{\mp}+\hbar}\left(k_2^{\mp}(v)
+f^{\mp}(v)k_1^{\mp}(v)e^{\mp}(v)\right)k_1^{\pm}(u)=0 
\end{eqnarray}

then using (\ref{kf1})(\ref{kf2}) and the following relations which also from 
(\ref{rll1})(\ref{rll2})
 and unitarity of $R$-matirx:
\begin{eqnarray}
&&k_1^{\pm}(u)e^{\pm}(u)k_1^{\pm}(v)-\frac{u-v}{u-v+\hbar}k_1^{\pm}(v)
k_1^{\pm}(u)e^{\pm}(u) \nonumber \\
&&{\hspace{4cm}}-\frac{\hbar}{u-v+\hbar}k_1^{\pm}(v)e^{\pm}(v)
k_1^{\pm}(u)=0  \\
&&k_1^{\pm}(u)e^{\pm}(u)k_1^{\mp}(v)-\frac{u_{\pm}-v_{\mp}}{u_{\pm}-v_{\mp}
+\hbar}k_1^{\mp}(v)k_1^{\pm}(u)e^{\pm}(u)  \nonumber \\
&&{\hspace{4cm}}-\frac{\hbar}{u_{\pm}-v_{\mp}
+\hbar}k_1^{\mp}(v)e^{\mp}(v)k_1^{\pm}(u)=0 
\end{eqnarray}
we obtain:
\begin{eqnarray}
&&e^{\pm}(u)f^{\pm}(v)+f^{\pm}(v)e^{\pm}(u)=\frac{\hbar}{u-v}k_1^{\pm}
(u)^{-1}k_2^{\pm}(u)-\frac{\hbar}{u-v}k_1^{\pm}(v)^{-1}k_2^{\pm}(v) \\
&&e^{\pm}(u)f^{\mp}(v)+f^{\mp}(v)e^{\pm}(u)=\frac{\hbar}{u_{\mp}-v_{\pm}}
k_1^{\pm}(u)^{-1}k_2^{\pm}(u)-\frac{\hbar}{u_{\pm}-v_{\mp}}k_1^{\mp}(v)^{-1}
k_2^{\mp}(v)
\end{eqnarray}
then 
\begin{eqnarray}
&&\{X^+(u) ~,~ X^-(v)\}=\hbar\left[ \delta(u_--v_+)k_1^+(u_-)^{-1}k_2^+(u_-)
-\delta(u_+-v_-)k_1^-(v_-)^{-1}k_2^-(v_-)\right]
\end{eqnarray}
where $ \delta(u-v)=\sum _{k \in Z}u^{k} v^{-k-1} $. 

In a word, we get all relations between $k_i^{\pm}(u)$ (i=1,2) and $ X^{\pm}(v)$:
\begin{eqnarray}
&&[k_1^{\pm}(u) ~,~k_1^{\pm}(v)]=[k_1^+(u) ~,~k_1^-(v)]=0  \\
&&[k_1^{\pm}(u) ~,~k_2^{\pm}(v)]=[k_2^{\pm}(u) ~,~k_2^{\pm}(v)]=0 \\
&&\frac{u_{\pm}-v_{\mp}}{u_{\pm}-v_{\mp}+\hbar}k_1^{\pm}(u)k_2^{\mp}(v)^{-1}
=\frac{u_{\mp}-v_{\pm}}{u_{\mp}-v_{\pm}+\hbar}k_2^{\mp}(v)^{-1}k_1^{\pm}(u) \\
&&\frac{u_--v_+-\hbar}{u_--v_++\hbar}k_2^+(u)^{-1}k_2^-(v)^{-1}
=\frac{u_+-v_--\hbar}{u_+-v_-+\hbar}k_2^-(v)^{-1}k_2^+(u)^{-1}  \\
&&k_i^{\pm}(u)^{-1}X^+(v)k_i^{\pm}(u)=\frac{u_{\pm}-v+\hbar}{u_{\pm}-v}X^+(v)
\hspace{1cm} (i=1,2) \\
&&k_i^{\pm}(u)X^-(v)k_i^{\pm}(u)^{-1}=\frac{u_{\mp}-v+\hbar}{u_{\mp}-v}X^-(v)
\hspace{1cm} (i=1,2) \\
&&\{X^+(u) ~,~X^+(v)\}=\{X^-(u) ~,~X^-(v)\}=0   \\
&&\{X^+(u) ~,~ X^-(v)\}=\hbar\left[ \delta(u_--v_+)k_1^+(u_-)^{-1}k_2^+(u_-)
-\delta(u_+-v_-)k_1^-(v_-)^{-1}k_2^-(v_-)\right]
\end{eqnarray}

Introducing a transformation for currents:
\begin{eqnarray}
&K^{\pm}(u)=k_1^{\pm}(u+\frac{\hbar}{2})^{-1}k_2^{\pm}(u+\frac{\hbar}{2})
&\hspace{1cm} H^{\pm}(u)=k_1^{\pm}(u-\frac{\hbar}{2})k_2^{\pm}
(u+\frac{\hbar}{2})  \\
&E(u)=X^+(u+\frac{\hbar}{2})&\hspace{1cm} F(u)=X^-(u+\frac{\hbar}{2})
\end{eqnarray}
 then all relations  among $K^{\pm}(u)$, $H^{\pm}(u)$, $E(u)$ and $F(u)$
can be calculated as follows:

\begin{eqnarray}
&&[K^{\pm}(u) ~,~ K^{\pm}(v)]=[K^+(u) ~,~K^-(v)]=0  \\
&&[H^{\pm}(u) ~,~ H^{\pm}(v)]=[K^{\pm}(u) ~,~ H^{\pm}(v)]=0  \\
&&\frac{u_{\mp}-v_{\pm}-\hbar}{u_{\mp}-v_{\pm}+\hbar}H^{\pm}(u)K^{\mp}(v)
=K^{\mp}(v)H^{\pm}(u)\frac{u_{\pm}-v_{\mp}-\hbar}{u_{\pm}-v_{\mp}+\hbar} \\
&&\left(\frac{u_{\mp}-v_{\pm}-\hbar}{u_{\mp}-v_{\pm}+\hbar}\right)^2
H^{\pm}(u)H^{\mp}(v)
=H^{\mp}(v)H^{\pm}(u)\left(\frac{u_{\pm}-v_{\mp}-\hbar}{u_{\pm}-v_{\mp}+\hbar}
\right)^2 \\
&&[K^{\pm}(u) ~,~E(v)]=[K^{\pm}(u) ~,~F(v)]=0  \\
&&E(v)H^{\pm}(u)=\frac{u_{\pm}-v+\hbar}{u_{\pm}-v-\hbar}H^{\pm}(u)E(v)  \\
&&H^{\pm}(u)F(v)=\frac{u_{\mp}-v+\hbar}{u_{\mp}-v-\hbar}F(v)H^{\pm}(u)  \\
&&\{E(u) ~,~ F(v)\}=\hbar\left[ \delta(u_--v_+)K^+(u_-)
-\delta(u_+-v_-)K^-(v_-)\right]
\end{eqnarray}

The co-product structure can be derived directly from (\ref{cop})(\ref{es}) :

\begin{eqnarray}
&&\triangle \left(E^{\pm}(u)\right)= E^{\pm}(u)\otimes 1+ 
H^{\pm}(u \mp \frac{\hbar}{4}c_1)\otimes E^{\pm}(u\mp \frac{\hbar}{2}c_1) \\
&&\triangle \left(F^{\pm}(u)\right)=1\otimes F^{\pm}(u)+
F^{\pm}(u \pm \frac{\hbar}{2} c_2)
\otimes H^{\pm}(u  \pm \frac{\hbar}{4}c_2) \\
&&\triangle \left(H^{\pm}(u)\right)= H^{\pm}(u\pm \frac{\hbar}{4}c_2)
\otimes H^{\pm}(u\mp \frac{\hbar}{4}c_1)-2F^{\pm}(u \pm \frac{\hbar}{4}c_2
\mp \frac{\hbar}{4}c_1-\hbar) H^{\pm}(u \pm \frac{\hbar}{4}c_2)\otimes 
 \nonumber \\
&&{\hspace{3cm}}H(u\mp \frac{\hbar}{4}c_1)^{\pm}E^{\pm}
(u \pm \frac{\hbar}{4}c_2 \mp \frac{\hbar}{4}c_1-\hbar)) \\
&&\triangle \left(K^{\pm}(u)\right)=K^{\pm}(u \pm \frac{\hbar}{4}c_2)\otimes 
K^{\pm}(u\mp \frac{\hbar}{4}c_1) \\
&&\epsilon (K^{\pm}(u)) =\epsilon (H^{\pm}(u)) =1 \\
&&\epsilon (E^{\pm}(u)) =\epsilon (F^{\pm}(u)) =0 
\end{eqnarray}

The quantum affine superalgebra $U_q\widehat{(gl(1|1))}$ corresponding to 
the trigonometric solution of super YBE have also been studied in our 
paper \cite{C2} and the  general
 case $U_q\widehat{(gl(m|n))}$ have been independently worked out
  by Zhang  \cite{YZZ}.

\end{document}